\documentclass[11pt]{amsart}
\usepackage{amssymb,mathrsfs}
\title{Problems in the Geometry of the Siegel-Jacobi Space}

\begin{document}

\author{Jae-Hyun Yang}

\address{Yang Institute for Advanced Study
\newline\indent
Hyundai 41 Tower, No. 1905
\newline\indent
293 Mokdongdong-ro, Yangcheon-gu
\newline\indent
Seoul 07997, Korea}
\address{Department of Mathematics
\newline\indent
Inha University
\newline\indent
Incheon 22212, Korea}

\email{jhyang@inha.ac.kr\ \ or\ \ jhyang8357@gmail.com}

\newtheorem{theorem}{Theorem}[section]
\newtheorem{lemma}{Lemma}[section]
\newtheorem{proposition}{Proposition}[section]
\newtheorem{remark}{Remark}[section]
\newtheorem{definition}{Definition}[section]

\renewcommand{\theequation}{\thesection.\arabic{equation}}
\renewcommand{\thetheorem}{\thesection.\arabic{theorem}}
\renewcommand{\thelemma}{\thesection.\arabic{lemma}}
\newcommand{\BR}{\mathbb R}
\newcommand{\BQ}{\mathbb Q}
\newcommand{\BT}{\mathbb T}
\newcommand{\BM}{\mathbb M}
\newcommand{\bn}{\bf n}
\def\charf {\mbox{{\text 1}\kern-.24em {\text l}}}
\newcommand{\BC}{\mathbb C}
\newcommand{\BZ}{\mathbb Z}

\thanks{\noindent{2010 Mathematics Subject Classification:} 14G35, 32F45, 32M10, 32Wxx, 53C22\\
\indent Keywords and phrases: Siegel-Jacobi space, Invariant metrics, Laplace operator,
Invariant differential \\
\indent operators, compactification.}

\begin{abstract}
{The Siegel-Jacobi space is a non-symmetric homogeneous space which is very important geometrically and arithmetically.
In this short paper, we propose the basic problems in the geometry of the Siegel-Jacobi space.}
\end{abstract}

\maketitle

\newcommand\tr{\triangleright}
\newcommand\al{\alpha}
\newcommand\be{\beta}
\newcommand\g{\gamma}
\newcommand\gh{\Cal G^J}
\newcommand\G{\Gamma}
\newcommand\de{\delta}
\newcommand\e{\epsilon}
\newcommand\z{\zeta}
\newcommand\vth{\vartheta}
\newcommand\vp{\varphi}
\newcommand\om{\omega}
\newcommand\p{\pi}
\newcommand\la{\lambda}
\newcommand\lb{\lbrace}
\newcommand\lk{\lbrack}
\newcommand\rb{\rbrace}
\newcommand\rk{\rbrack}
\newcommand\s{\sigma}
\newcommand\w{\wedge}
\newcommand\fgj{{\frak g}^J}
\newcommand\lrt{\longrightarrow}
\newcommand\lmt{\longmapsto}
\newcommand\lmk{(\lambda,\mu,\kappa)}
\newcommand\Om{\Omega}
\newcommand\ka{\kappa}
\newcommand\ba{\backslash}
\newcommand\ph{\phi}
\newcommand\M{{\Cal M}}
\newcommand\bA{\bold A}
\newcommand\bH{\bold H}
\newcommand\D{\Delta}

\newcommand\Hom{\text{Hom}}
\newcommand\cP{\Cal P}

\newcommand\cH{\Cal H}

\newcommand\pa{\partial}

\newcommand\pis{\pi i \sigma}
\newcommand\sd{\,\,{\vartriangleright}\kern -1.0ex{<}\,}
\newcommand\wt{\widetilde}
\newcommand\fg{\frak g}
\newcommand\fk{\frak k}
\newcommand\fp{\frak p}
\newcommand\fs{\frak s}
\newcommand\fh{\frak h}
\newcommand\Cal{\mathcal}

\newcommand\fn{{\frak n}}
\newcommand\fa{{\frak a}}
\newcommand\fm{{\frak m}}
\newcommand\fq{{\frak q}}
\newcommand\CP{{\mathcal P}_n}
\newcommand\Hnm{{\mathbb H}_n \times {\mathbb C}^{(m,n)}}
\newcommand\BD{\mathbb D}
\newcommand\BH{\mathbb H}
\newcommand\CCF{{\mathcal F}_n}
\newcommand\CM{{\mathcal M}}
\newcommand\Gnm{\Gamma_{n,m}}
\newcommand\Cmn{{\mathbb C}^{(m,n)}}
\newcommand\Yd{{{\partial}\over {\partial Y}}}
\newcommand\Vd{{{\partial}\over {\partial V}}}

\newcommand\Ys{Y^{\ast}}
\newcommand\Vs{V^{\ast}}
\newcommand\LO{L_{\Omega}}
\newcommand\fac{{\frak a}_{\mathbb C}^{\ast}}

\vskip 7mm

\centerline{\large \bf Table of Contents}

\vskip 0.5cm $ \qquad\qquad\qquad\qquad\textsf{\large \ 1.
Introduction}$\vskip 0.021cm

$\qquad\qquad\qquad\qquad \textsf{\large\ 2. Brief Review on the Geometry of the Siegel Space }$
\vskip 0.021cm

$ \qquad\qquad\qquad\qquad  \textsf{\large\ 3. Basic Problems in the Geometry of the Siegel-Jacobi Space }$
\vskip 0.021cm

$ \qquad\qquad\qquad\qquad \textsf{\large\ 4. Final Remarks}$
\vskip 0.021cm


$ \qquad\qquad\qquad\qquad\textsf{\large References }$

\vskip 10mm

\begin{section}{{\bf Introduction}}
\setcounter{equation}{0} For a given fixed positive integer $n$,
we let
$${\mathbb H}_n=\,\{\,\Omega\in \BC^{(n,n)}\,|\ \Om=\,^t\Om,\ \ \ \text{Im}\,\Om>0\,\}$$
be the Siegel upper half plane of degree $n$ and let
$$Sp(n,\BR)=\{ M\in \BR^{(2n,2n)}\ \vert \ ^t\!MJ_nM= J_n\ \}$$
be the symplectic group of degree $n$, where $F^{(k,l)}$ denotes
the set of all $k\times l$ matrices with entries in a commutative
ring $F$ for two positive integers $k$ and $l$, $^t\!M$ denotes
the transposed matrix of a matrix $M$ and
$$J_n=\begin{pmatrix} 0&I_n\\
                   -I_n&0\end{pmatrix}.$$
Then $Sp(n,\BR)$ acts on $\BH_n$ transitively by
\begin{equation}
M\cdot\Om=(A\Om+B)(C\Om+D)^{-1},
\end{equation} where $M=\begin{pmatrix} A&B\\
C&D\end{pmatrix}\in Sp(n,\BR)$ and $\Om\in \BH_n.$ Let
$$\G_n=Sp(n,\BZ)=\left\{ \begin{pmatrix} A&B\\
C&D\end{pmatrix}\in Sp(n,\BR) \,\big| \ A,B,C,D\
\textrm{integral}\ \right\}$$ be the Siegel modular group of
degree $n$. This group acts on $\BH_n$ properly discontinuously.
C. L. Siegel investigated the geometry of $\BH_n$ and automorphic
forms on $\BH_n$ systematically. Siegel\,\cite{Si1} found a
fundamental domain ${\mathcal F}_n$ for $\G_n\ba\BH_n$ and
described it explicitly. Moreover he calculated the volume of
$\CCF.$ We also refer to \cite{M2},\,\cite{Si1} for
some details on $\CCF.$

\vskip 0.1cm For two
positive integers $m$ and $n$, we consider the Heisenberg group
$$H_{\BR}^{(n,m)}=\big\{\,(\la,\mu;\ka)\,|\ \la,\mu\in \BR^{(m,n)},\ \kappa\in\BR^{(m,m)},\
\ka+\mu\,^t\la\ \text{symmetric}\ \big\}$$ endowed with the
following multiplication law
$$\big(\la,\mu;\ka\big)\circ \big(\la',\mu';\ka'\big)=\big(\la+\la',\mu+\mu';\ka+\ka'+\la\,^t\mu'-
\mu\,^t\la'\big)$$
with $\big(\la,\mu;\ka\big),\big(\la',\mu';\ka'\big)\in H_{\BR}^{(n,m)}.$
We define the {\it Jacobi group} $G^J$ of degree $n$ and index $m$ that is the semidirect product of
$Sp(n,\BR)$ and $H_{\BR}^{(n,m)}$
$$G^J=Sp(n,\BR)\ltimes H_{\BR}^{(n,m)}$$
endowed with the following multiplication law
$$
\big(M,(\lambda,\mu;\kappa)\big)\cdot\big(M',(\lambda',\mu';\kappa'\,)\big)
=\, \big(MM',(\tilde{\lambda}+\lambda',\tilde{\mu}+ \mu';
\kappa+\kappa'+\tilde{\lambda}\,^t\!\mu'
-\tilde{\mu}\,^t\!\lambda'\,)\big)$$ with $M,M'\in Sp(n,\BR),
(\lambda,\mu;\kappa),\,(\lambda',\mu';\kappa') \in
H_{\BR}^{(n,m)}$ and
$(\tilde{\lambda},\tilde{\mu})=(\lambda,\mu)M'$. Then $G^J$ acts
on $\BH_n\times \BC^{(m,n)}$ transitively by
\begin{equation}
\big(M,(\lambda,\mu;\kappa)\big)\cdot
(\Om,Z)=\Big(M\cdot\Om,(Z+\lambda \Om+\mu)
(C\Omega+D)^{-1}\Big), \end{equation} where $M=\begin{pmatrix} A&B\\
C&D\end{pmatrix} \in Sp(n,\BR),\ (\lambda,\mu; \kappa)\in
H_{\BR}^{(n,m)}$ and $(\Om,Z)\in \BH_n\times \BC^{(m,n)}.$ We note
that the Jacobi group $G^J$ is {\it not} a reductive Lie group and
the homogeneous space ${\mathbb H}_n\times \BC^{(m,n)}$ is not a
symmetric space. From now on, for brevity we write
$\BH_{n,m}=\BH_n\times \BC^{(m,n)}.$ The homogeneous space $\BH_{n,m}$ is called the
{\it Siegel-Jacobi space} of degree $n$ and index $m$.

\vskip 0.21cm
In this short article, we propose the basic and natural problems in the geometry
of the Siegel-Jacobi space.

\vskip 0.51cm \noindent {\bf Notations:} \ \ We denote by
$\BQ,\,\BR$ and $\BC$ the field of rational numbers, the field of
real numbers and the field of complex numbers respectively. We
denote by $\BZ$ the ring of integers. The symbol ``:='' means that
the expression on the right is the definition of that on the left.
For two positive integers $k$ and $l$, $F^{(k,l)}$ denotes the set
of all $k\times l$ matrices with entries in a commutative ring
$F$. For a square matrix $A\in F^{(k,k)}$ of degree $k$,
$\sigma(A)$ denotes the trace of $A$. For any $M\in F^{(k,l)},\
^t\!M$ denotes the transpose of a matrix $M$. $I_n$ denotes the
identity matrix of degree $n$. For a complex matrix $A$,
${\overline A}$ denotes the complex {\it conjugate} of $A$.
For a number field $F$, we denote
by ${\mathbb A}_F$ the ring of adeles of $F$. If $F=\BQ$, the
subscript will be omitted.

\vskip 0.1cm

\end{section}

\begin{section}{{\bf Brief Review on the Geometry of the Siegel Space}}
\setcounter{equation}{0}

\newcommand\POB{ {{\partial}\over {\partial{\overline \Omega}}} }
\newcommand\PZB{ {{\partial}\over {\partial{\overline Z}}} }
\newcommand\PX{ {{\partial}\over{\partial X}} }
\newcommand\PY{ {{\partial}\over {\partial Y}} }
\newcommand\PU{ {{\partial}\over{\partial U}} }
\newcommand\PV{ {{\partial}\over{\partial V}} }
\newcommand\PO{ {{\partial}\over{\partial \Omega}} }
\newcommand\PZ{ {{\partial}\over{\partial Z}} }

\vskip 0.21cm
We let $G:=Sp(n,\BR)$ and $K=U(n).$ The stabilizer of the action (1.1) at $iI_n$
is
\begin{equation*}
  \left\{ \begin{pmatrix} \,A & B \\ -B & A \end{pmatrix} \Big| \ A+iB\in U(n)\,\right\}
  \cong U(n).
\end{equation*}
Thus we get the biholomorphic map
\begin{equation*}
G/K \lrt \BH_n, \qquad gK \mapsto g\!\cdot\! iI_n,  \quad g\in G.
\end{equation*}
$\BH_n$ is a Hermitian symmetric manifold.

\vskip 0.21cm For $\Om=(\omega_{ij})\in\BH_n,$ we write $\Om=X+iY$
with $X=(x_{ij}),\ Y=(y_{ij})$ real. We put $d\Om=(d\om_{ij})$ and $d{\overline\Om}=(d{\overline\om}_{ij})$. We
also put
$$\PO=\,\left(\,
{ {1+\delta_{ij}}\over 2}\, { {\!\!\partial}\over {\partial \om_{ij} }
} \,\right) \qquad\text{and}\qquad \POB=\,\left(\, {
{1+\delta_{ij}}\over 2}\, { {\!\!\partial}\over {\partial {\overline
{\om}}_{ij} } } \,\right).$$ C. L. Siegel \cite{Si1} introduced
the symplectic metric $ds_{n;A}^2$ on $\BH_n$ invariant under the action
(1.1) of $Sp(n,\BR)$ that is given by
\begin{equation}
ds_{n;A}^2=A\,\s (Y^{-1}d\Om\, Y^{-1}d{\overline\Om}),\qquad A>0.
\end{equation}
It is known that the metric $ds_{n;A}^2$ is a K{\"a}hler-Einstein metric.
H. Maass \cite{M1} proved that its Laplace operator $\Delta_{n;A}$ is given by
\begin{equation}
\Delta_{n;A}=\,{4\over A}\,\s \left(\,Y\,
{}^{{}^{{}^{{}^\text{\scriptsize $t$}}}}\!\!\!
\left(Y\POB\right)\PO\right).\end{equation} And
\begin{equation}
dv_n(\Om)=(\det Y)^{-(n+1)}\prod_{1\leq i\leq j\leq n}dx_{ij}\,
\prod_{1\leq i\leq j\leq n}dy_{ij}\end{equation} is a
$Sp(n,\BR)$-invariant volume element on
$\BH_n$\,(cf.\,\cite[p.\,130]{Si2}).

\vskip 2mm
Siegel proved the following theorem for the Siegel space $(\BH_n, ds^2_{n;1}).$
\begin{theorem}\,({\bf Siegel\,\cite{Si1}}).
(1) There exists exactly one geodesic joining two arbitrary points
$\Om_0,\,\Om_1$ in $\BH_n$. Let $R(\Om_0,\Om_1)$ be the
cross-ratio defined by
\begin{equation*}
R(\Om_0,\Om_1)=(\Om_0-\Om_1)(\Om_0-{\overline
\Om}_1)^{-1}(\overline{\Om}_0-\overline{\Om}_1)(\overline{\Om}_0-\Om_1)^{-1}.
\end{equation*}
For brevity, we put $R_*=R(\Om_0,\Om_1).$ Then the symplectic
length $\rho(\Om_0,\Om_1)$ of the geodesic joining $\Om_0$ and
$\Om_1$ is given by
\begin{equation*}
\rho(\Om_0,\Om_1)^2=\s \left( \left( \log { {1+R_*^{\frac 12}
}\over {1-R_*^{\frac 12} } }\right)^2\right),
\end{equation*} where
\begin{equation*}
\left( \log { {1+R_*^{\frac 12} }\over {1-R_*^{\frac 12} }
}\right)^2=\,4\,R_* \left( \sum_{k=0}^{\infty} { {R_*^k}\over
{2k+1}}\right)^2.
\end{equation*}

\noindent (2) For $M\in Sp(n,\BR)$, we set
$${\tilde \Om}_0=M\cdot \Om_0\quad \textrm{and}\quad {\tilde \Om}_1=M\cdot
\Om_1.$$ Then $R(\Om_1,\Om_0)$ and
$R({\tilde\Om}_1,{\tilde\Om}_0)$ have the same eigenvalues.

\vskip 2mm\noindent
\noindent (3) All geodesics are symplectic images of the special
geodesics
\begin{equation*}
\alpha(t)=i\,\textrm{diag}(a_1^t,a_2^t,\cdots,a_n^t),
\end{equation*}
where $a_1,a_2,\cdots,a_n$ are arbitrary positive real numbers
satisfying the condition
$$\sum_{k=1}^n \left( \log a_k\right)^2=1.$$
\end{theorem}
\noindent The proof of the above theorem can be found in
\cite[pp.\,289-293]{Si1}.

\vskip 5mm
Let $\mathbb D(\BH_n)$ be the algebra of all differential operators on $\BH_n$ invariant
under the action (1.1). Then according to Harish-Chandra
\cite{HC1, HC2},
$$  \mathbb D (\BH_n)= \BC [ D_1,\cdots,D_n ],$$
where $D_1,\cdots,D_n$ are algebraically independent invariant differential operators on $\BH_n$.
That is, $\mathbb D (\BH_n)$ is a commutative algebra that is finitely generated by $n$
algebraically independent invariant differential operators on $\BH_n$. Maass \cite{M2} found the explicit
$D_1,\cdots,D_n$. Let $\frak g_\BC$ be the complexification of the Lie algebra of $G$.
It is known that $\BD(\BH_n)$ is isomorphic to the
center of the universal enveloping algebra of $\fg_{\BC}$\,(cf.\,\cite{He}).

\vskip 5mm\noindent
{\bf Example.} We consider the simplest case $n=1$ and $A=1.$  Let $\BH$ be the Poincar{\'e} upper half plane. Let $\omega=x+iy\in \BH$ with $x,y\in \BR$ and $y>0.$ Then the Poincar{\'e} metric
$$ ds^2= {{dx^2 +dy^2}\over {y^2} } = {{d\omega\, d{\overline \omega} }\over {y^2} }$$
is a $SL(2,\BR)$-invariant K{\"a}hler-Einstein metric on $\BH$.
The geodesics of $(\BH, ds^2)$ are either straight vertical lines perpendicular to the $x$-axis or
circular arcs perpendicular to the $x$-axis (half-circles whose origin is on the $x$-axis).
The Laplace operator $\Delta$ of $(\BH,ds^2)$ is given by
$$\Delta=y^2 \left( {{\partial^2\,\, }\over {\partial x^2}} + {{\partial^2\,\, }\over {\partial y^2}}
\right)$$
and
$$ dv= {{dx\wedge dy}\over {y^2}}$$
is a $SL(2,\BR)$-invariant volume element. The scalar curvature, i.e., the Gaussian curvature is $-1.$
The algebra $\BD(\BH)$ of all $SL(2,\BR)$-invariant differential operators on $\BH$ is given by
$$\BD(\BH)=\BC [\Delta].$$

The distance between two points $\omega_1=x_1+iy_1$ and $\omega_2=x_2+iy_2$ in $(\BH, ds^2)$ is given by
\begin{eqnarray*}
  \rho (\omega_1,\omega_2) &=& 2\,{\rm {ln}}
  { { \sqrt{(x_2-x_1)^2 + (y_2-y_1)^2}+\sqrt{(x_2-x_1)^2 + (y_2+y_1)^2} } \over  {2 \sqrt{y_1y_2} } } \\
  &=& \cosh^{-1}\left( 1 + { {(x_2-x_1)^2 + (y_2-y_1)^2}\over {2y_1y_2} } \right) \\
   &=& 2\, \sinh^{-1} {\frac 12} \sqrt{ { {(x_2-x_1)^2 + (y_2-y_1)^2}\over {y_1y_2} } }.
\end{eqnarray*}

\end{section}

\begin{section}{{\bf Basic Problems in the Geometry of the Siegel-Jacobi Space}}
\setcounter{equation}{0}

\newcommand\POB{ {{\partial}\over {\partial{\overline \Omega}}} }
\newcommand\PZB{ {{\partial}\over {\partial{\overline Z}}} }
\newcommand\PX{ {{\partial}\over{\partial X}} }
\newcommand\PY{ {{\partial}\over {\partial Y}} }
\newcommand\PU{ {{\partial}\over{\partial U}} }
\newcommand\PV{ {{\partial}\over{\partial V}} }
\newcommand\PO{ {{\partial}\over{\partial \Omega}} }
\newcommand\PZ{ {{\partial}\over{\partial Z}} }

\vskip 0.2cm
For a coordinate
$(\Om,Z)\in \BH_{n,m}$ with $\Om=(\omega_{\mu\nu})$ and
$Z=(z_{kl})$, we put $d\Om,\,d{\overline \Om},\,\PO,\,\POB$ as
before and set
\begin{eqnarray*}
Z\,&=&U\,+\,iV,\quad\ \ U\,=\,(u_{kl}),\quad\ \ V\,=\,(v_{kl})\ \
\text{real},\\
dZ\,&=&\,(dz_{kl}),\quad\ \ d{\overline Z}=(d{\overline z}_{kl}),
\end{eqnarray*}
$$\PZ=\begin{pmatrix} {\partial}\over{\partial z_{11}} & \hdots &
 {\partial}\over{\partial z_{m1}} \\
\vdots&\ddots&\vdots\\
 {\partial}\over{\partial z_{1n}} &\hdots & {\partial}\over
{\partial z_{mn}} \end{pmatrix},\quad \PZB=\begin{pmatrix}
{\partial}\over{\partial {\overline z}_{11} }   &
\hdots&{ {\partial}\over{\partial {\overline z}_{m1} }  }\\
\vdots&\ddots&\vdots\\
{ {\partial}\over{\partial{\overline z}_{1n} }  }&\hdots &
 {\partial}\over{\partial{\overline z}_{mn} }  \end{pmatrix}.$$

\newcommand\bw{d{\overline W}}
\newcommand\bz{d{\overline Z}}
\newcommand\bo{d{\overline \Omega}}

\vskip 0.3cm
 The author proved the following theorems in \cite{Y1}.

\begin{theorem} For any two positive real numbers
$A$ and $B$,
\begin{eqnarray}
ds_{n,m;A,B}^2&=&\,A\, \sigma\Big( Y^{-1}d\Om\,Y^{-1}d{\overline
\Om}\Big) \nonumber \\
&& \ \ + \,B\,\bigg\{ \sigma\Big(
Y^{-1}\,^tV\,V\,Y^{-1}d\Om\,Y^{-1} d{\overline \Om} \Big)
 +\,\sigma\Big( Y^{-1}\,^t(dZ)\,\bz\Big) \nonumber\\
&&\quad\quad -\sigma\Big( V\,Y^{-1}d\Om\,Y^{-1}\,^t(\bz)\Big)\,
-\,\sigma\Big( V\,Y^{-1}d{\overline \Om}\, Y^{-1}\,^t(dZ)\,\Big)
\bigg\} \nonumber
\end{eqnarray}
is a Riemannian metric on $\BH_{n,m}$ which is invariant under the action (1.2) of $G^J.$
\end{theorem}
\vskip 1mm\noindent
{\it Proof.} See Theorem 1.1 in \cite{Y1}. \hfill $\Box$

\vskip 3mm

\begin{theorem} The Laplace operator $\Delta_{m,m;A,B}$ of the $G^J$-invariant metric $ds_{n,m;A,B}^2$ is given by
\begin{equation}
\Delta_{n,m;A,B}=\,{\frac 4A}\,{\mathbb M}_1 + {\frac 4B}
{\mathbb M}_2,
\end{equation}
where
\begin{eqnarray*}
{\mathbb M}_1\,&=&  \sigma\left(\,Y\,
{}^{{}^{{}^{{}^\text{\scriptsize $t$}}}}\!\!\!
\left(Y\POB\right)\PO\,\right)\, +\,\sigma\left(\,VY^{-1}\,^tV\,
{}^{{}^{{}^{{}^\text{\scriptsize $t$}}}}\!\!\!
\left(Y\PZB\right)\,\PZ\,\right)\\
& &\ \
+\,\sigma\left(V\,
{}^{{}^{{}^{{}^\text{\scriptsize $t$}}}}\!\!\!
\left(Y\POB\right)\PZ\,\right)
+\,\sigma\left(\,^tV\,
{}^{{}^{{}^{{}^\text{\scriptsize $t$}}}}\!\!\!
\left(Y\PZB\right)\PO\,\right)\nonumber
\end{eqnarray*}
and
\begin{equation*}
{\mathbb M}_2=\,\sigma\left(\, Y\,\PZ\,
{}^{{}^{{}^{{}^\text{\scriptsize $t$}}}}\!\!\!
\left(
\PZB\right)\,\right).
\end{equation*}
Furthermore ${\mathbb M}_1$ and ${\mathbb M}_2$ are differential operators on $\BH_{n,m}$ invariant under the action (1.2) of $G^J.$
\end{theorem}
\vskip 1mm\noindent
{\it Proof.} See Theorem 1.2 in \cite{Y1}. \hfill $\Box$

\vskip 3mm
\begin{remark}
Erik Balslev \cite{B} developed the spectral theory of $\Delta_{1,1;1,1}$ on $\Bbb H_{1,1}$ for certain arithmetic subgroups of the Jacobi modular group to
prove that the set of all eigenvalues of $\Delta_{1,1;1,1}$ satisfies the Weyl law.
\end{remark}

\begin{remark}
The scalar curvature of $(\BH_{1,1}, ds^2_{1,1;A,B})$ is $-{3\over A}$ and hence is independent of the parameter $B$. We refer to \cite{Y3}
for more detail.
\end{remark}

\begin{remark}
Yang and Yin \cite{YY} showed that $ds^2_{n,m;A,B}$ is a K{\"a}hler metric.
For some applications of the invariant metric $ds^2_{n,m;A,B}$ we refer to \cite{YY}.
\end{remark}

\vskip 5mm
Now we propose the basic and natural problems.

\vskip 3mm\noindent
{\bf Problem 1.} Find all the geodesics of $(\BH_{n,m},ds^2_{n,m;A,B})$ explicitly.

\vskip 3mm\noindent
{\bf Problem 2.} Compute the distance between two points $(\Omega_1,Z_1)$ and $(\Omega_2,Z_2)$ of $\BH_{n,m}$ explicitly.

\vskip 3mm\noindent
{\bf Problem 3.} Compute the Ricci curvature tensor and the scalar curvature of $(\BH_n,ds^2_{n,m;A,B})$.

\vskip 3mm\noindent
{\bf Problem 4.} Find all the eigenfunctions of the Laplace operator $\Delta_{n,m;A,B}$.

\vskip 3mm\noindent
{\bf Problem 5.} Develop the spectral theory of $\Delta_{n,m;A,B}$.

\vskip 3mm\noindent
{\bf Problem 6.} Describe the algebra of all $G^J$-invariant differential operators on $\BH_{n,m}$ explicitly. We refer to \cite{Y1+, Y2, YY} for some details.

\vskip 3mm\noindent
{\bf Problem 7.} Find the trace formula for the Jacobi group $G^J (\mathbb A)$.

\vskip 3mm\noindent
{\bf Problem 8.} Discuss the behaviour of the analytic torsion of the Siegel-Jacobi space $\BH_{n,m}$ or
the arithmetic quotients of $\BH_{n,m}$.

\vskip 5mm
We make some remarks on the above problems.
\begin{remark}
Problem 1 reduces to trying to solve a system of ordinary differential equations explicitly.
If Problem 2 is solved, the distance formula would be a very beautiful one that generalizes the distance
formula $\rho(\Omega_0,\Omega_1)$ given by Theorem 2.1 (the Siegel space case).
\end{remark}

\begin{remark}
Problem 3 was recently solved in the case that $n=1$ and $m$ is arbitrary. Precisely the
scalar and Ricci curvatures of the Siegel-Jacobi space $(\BH_{1,m}, ds^2_{1,m;A,B})\,(m\geq 1)$
were completely computed by G. Khan and J. Zhang \cite[Proposition 8, pp.\,825--826]{KZ}.
Furthermore Khan and Zhang proved that $(\BH_{1,m}, ds^2_{1,m;A,B})\,(m\geq 1)$ has  non-negative  orthogonal  anti-bisectional curvature (cf.\,\cite[Proposition 9, p.\,826]{KZ}).
\end{remark}

\begin{remark}
Concerning Problem 4 and Problem 5, computing eigenfunctions explicitly is a tall order, but if this can be done it will shed a lot of light onto the geometry of this space. And understanding the spectral geometry seems to be a central question which will likely have applications in number theory and other areas.
\end{remark}

\begin{remark}
The algebra $\BD(\BH_{n,m})$ of all $G^J$-invariant differential operators on $\BH_{n,m}$ is not commutative.
Concerning Problem 6, the case $n=m=1$ was completely solved by M. Itoh, H. Ochiai and J.-H. Yang in 2013.
They proved that the noncommutative algebra $\BD(\BH_{1,1})$ is generated by four explicit generators
$D_1,D_2,D_3,D_4$, and found the relations among those $D_i\,(1\leq i\leq 4).$
For more precise statements, we refer to \cite[pp.\,56--58]{Y1+} and \cite[pp.\,285--290]{Y2}.
We note that the above four generators $D_i\,(1\leq i\leq 4)$ are not algebraically independent.
\end{remark}

\begin{remark}
The solution of Problem 7 will provide lots of arithmetic properties of the Siegel-Jacobi space.
\end{remark}

\end{section}

\begin{section}{{\bf Final Remarks}}
\setcounter{equation}{0}

\vskip 2mm
Let $\Gamma_n (N)$ be the principal congruence subgroup of the Siegel modular group $\Gamma_n$.
Let $\mathfrak{X}_n (N):=\Gamma_n(N)\backslash \BH_n$ be the moduli of $n$-dimensional principally polarized abelian varieties with level $N$-structure. The Mumford school \cite{AMRT} found toroidal compactifications of $\mathfrak{X}_n (N)$ which are usefully applied in the study of the geometry and arithmetic of $\mathfrak{X}_n (N)$. D. Mumford \cite{Mf} proved the Hirzebruch's Proportionality Theorem in the non-compact case introducing a {\it good\ singular} Hermitian metric on an automorphic vector bundle on a smooth toroidal compactification of $\mathfrak{X}_n (N)$ with $N\geq 3.$
\vskip 3mm
We set
\begin{equation*}
  \Gamma_{n,m}(N):=\Gamma_n \ltimes H_{\BZ}^{(n,m)},
\end{equation*}
where
\begin{equation*}
   H_{\BZ}^{(n,m)}= \left\{ (\lambda,\mu;\kappa)\in H_{\BR}^{(n,m)}\,\vert \ \lambda,\mu,\kappa\
   integral\,\right\}.
\end{equation*}

Let
\begin{equation*}
  \mathfrak X_{n,m}(N):= \Gamma_{n,m}(N)\backslash \BH_{n,m}
\end{equation*}
be the universal abelian variety. An arithmetic toroidal compactification of $\mathfrak X_{n,m}(N)$ was intensively investigated by R. Pink \cite{P}. D. Mumford described very nicely a toroidal compactification of the universal elliptic curve $\mathfrak X_{1,1}(N)$ (cf.\,\cite[pp.\,14--25]{AMRT}). The geometry of $\mathfrak X_{n,m}(N)$ is closely related to the theory of Jacobi forms (cf. \cite{BKK, K1, K2, KvP}).
Jacobi forms play an important role in the study of the geometric and arithmetic of $\mathfrak X_{n,m}(N)$.
We refer to \cite{EZ, Z} for the theory of Jacobi forms.

\end{section}


\vspace{10mm}


\begin{thebibliography}{99}

\bibitem{AMRT} A. Ash, D. Mumford, M. Rapoport and Y.-S. Tai, {\em Smooth Compactifications of Locally Symmetric Varieties. With the Collaboration of Peter Scholze}, 2nd ed. Cambridge Math. Library, Cambridge Univ. Press, Cambridge (2010).

\bibitem{B} E. Balslev, {\em Spectral theory of the Laplacian on the modular Jacobi group manifold},
preprint, Aarhus University (2012).

\bibitem{BKK} J.\,I. Burgos Gil, J. Kramer and U. K{\"uhn}, {\em The Singularities of the invariant metric on the Jacobi line bundle}, in Recent Advances in Hodge Theory: Period Domains, Algebraic Cycles and
    Arithmetic, M. Kerr and G. Pearlstein (eds.), pp.\,45-77, London Mathematical Society, Lecture Note
    Series {\bf 427}, Cambridge University Press (February 2016).

\bibitem{EZ} M. Eichler and D. Zagier, {\em The Theory of Jacobi Forms}, Progress in Mathematics {\bf 55}, Birkh{\"a}user, Boston, Basel and Stuttgart, 1985.

\bibitem{HC1} Harish-Chandra, \textit{Representations of a semisimple Lie group on a Banach space. I.,} Trans. Amer. Math. Soc. {\bf 75} (1953), 185-243.

\bibitem{HC2} Harish-Chandra, \textit{The characters of semisimple Lie groups,} Trans. Amer. Math. Soc. {\bf 83} (1956), 98-163.

\bibitem{He} S. Helgason, {\em Groups and geometric analysis,} Academic Press, New York (1984).

\bibitem{KZ} G. Khan and J. Zhang, {\em A hall of statistical mirrors}, Asian J. Math. {\bf 26} (2022),
no.\,6, 809--846.

\bibitem{K1} J. Kramer, {\em A geometrical approach to the theory of Jacobi forms}, Compositio Math., {\bf 79}(1991), 1-19.

\bibitem{K2} J. Kramer, {\em A geometrical approach to Jacobi forms, revisited}, Oberwolfach Report 22/2014, 1270-1272.

\bibitem{KvP} J. Kramer and A.-M. von Pippich, {\em Snapshots of Modern Mathematics from Oberwohlfach:
Special Values of Zeta Functions and Areas of Triangles}, Notices of the AMS, vol.\,63, No.\,8 (September 2016), 917-922.



\bibitem{M1} H. Maass, {\em Die Differentialgleichungen in der Theorie der Siegelschen Modulfunktionen}, Math. Ann. {\bf 126} (1953), 44--68.

\bibitem{M2} H. Maass, {\em Siegel modular forms and Dirichlet series,} Lecture Notes in Math. {\bf 216}, Springer-Verlag, Berlin-Heidelberg-New York (1971).

\bibitem{Mf} D. Mumford, {\em Hirzebruch's Proportionality Theorem in the Non-Compact Case}, Invent. Math. {\bf 42} (1977), 239--272.

\bibitem{P} R. Pink, {\em Arithmetical compactification of mixed Shimura varieties}, Ph.\,D Thesis, Bonn 1989.

\bibitem{Si1} C.~L.~Siegel, \emph{Symplectic Geometry,} Amer. J. Math. {\bf 65} (1943), 1-86; Academic Press, New York and London (1964);
Gesammelte Abhandlungen, no.\,\,41, vol. II, Springer-Verlag (1966), 274-359.

\bibitem{Si2}  C. L. Siegel, {\em Topics in Complex Function Theory\,: Abelian Functions and Modular Functions of Several Variables}, vol. III, Wiley-Interscience, 1973.


\bibitem{Y1} J.-H. Yang, {\em Invariant metrics and Laplacians on Siegel-Jacobi space,} Journal of Number Theory {\bf 127} (2007), 83--102.

\bibitem{Y1+} J.-H. Yang, {\em Invariant differential operators on Siegel-Jacobi space and Maass-Jacobi forms}, Proceedings of the International Conference (October 10--12, 2012) on Geometry, number theory and representation theory, 37–63. KM Kyung Moon Sa, Seoul (2013).
ISBN:978-89-88615-35-5

\bibitem{Y2} J.-H. Yang, \emph{Geometry and Arithmetic on the Siegel-Jacobi Space}, Geometry and Analysis on Manifolds, In Memory of Professor Shoshichi Kobayashi (edited by T. Ochiai, A. Weinstein et al), Progress in Mathematics, Volume {\bf 308}, Birkh{\"a}user, Springer International Publishing AG Switzerland (2015), 275--325.

\bibitem{Y3} J.-H. Yang, Y.-H. Yong, S.-N. Huh, J.-H. Shin and G.-H. Min,
{\em Sectional Curvatures of the Siegel-Jacobi Space,} Bull. Korean Math. Soc. {\bf 50} (2013), No. 3, pp. 787-799.

\bibitem{YY} J. Yang and L. Yin, {\em Differetial operators for Siegel-Jacobi forms,} Science China Mathematics, Vol. {\bf 59} (2016), No.~6, pp. 1029--1050.

\bibitem{Z}  C. Ziegler, {\em Jacobi Forms of Higher Degree}, Abh. Math. Sem. Hamburg {\bf 59} (1989), 191--224.


\end{thebibliography}
\end{document}